\documentclass{article}
\usepackage{latexsym}

\title{\bf An Extension of Barbashin-Krasovski-LaSalle Theorem to a Class of 
 Nonautonomous Systems}
\author{Radu Balan\footnote{This is a late submission of an earlier report. 
Original version date: October 19, 1995, when the author was with Princeton University.} \\ 
 \small Siemens Corporate Research \\
 \small 755 College Road East \\
 \small Princeton, NJ 08540 \\
 \small e-mail: rvbalan@yahoo.com
 }

\newcommand{\R}{{\bf R}}

\newcommand{\bez}{{\bar{B}}_{{\varepsilon}_0}}
\newcommand{\ep}{\varepsilon}
\newcommand{\pp}{\parallel}

\newtheorem{Theorem}{THEOREM}

\newtheorem{Lemma}[Theorem]{LEMMA}
\newtheorem{Example}{EXAMPLE}

\begin{document}
\rm
\maketitle
\begin{abstract}

In this paper we  give an extension of the
 Barbashin-Krasovski-LaSalle Theorem 
 to a class of time-varying dynamical
 systems, namely the class of systems for which the restricted vector field to 
 the
 zero-set of the time derivative of the Liapunov function is time invariant
 and this set includes some trajectories.
  Our goal is to improve the sufficient conditions for the case of
 uniform asymptotic stability of the equilibrium.
 We obtain an extension of an well-known linear result
 to the case of zero-state detectability (given $(C,A)$ a detectable pair,
 if there exists a
 positive semidefinite matrix $P\geq 0$ such that: $A^TP+PA+C^TC=0$ then
 $A$ is Hurwitz -i.e. it has all eigenvalues with negative real part) as well
 as a result about robust stabilizability of nonlinear affine control
 systems.

\end{abstract}
\vspace{15mm}

{\bf \mbox{}~~Key words:}  Invariance Principle, Liapunov functions,
 detectability, robust stabilizability

\section{Introduction and Main Result}

Let us consider the following time-varying dynamical system:
\begin{equation}
\label{e1}
\dot{x}=f(t,x)~~~,~x\in D, t\in\R
\end{equation}
where $D$ is a domain in ${\R}^n$ containing the origin ($0\in
D\subset{\R}^n$). About $f$ we suppose the following:

1) $f(t,0)=0$, for any $t\in\R$;

2) Uniformly continuous in $t$, uniformly in $x\in D$, i.e. $\forall\ep>0
\exists {\delta}_{\ep}>0$ s.t. $\forall t_1,t_2\in\R,|t_1-t_2|<{\delta}_{\ep}$
 and $\forall x\in D,~\pp f(t_1,x)-f(t_2,x)\pp<\ep$;

3) Uniformly local Lipschitz continuous in $x$ for any $t\in\R$ , i.e.
 for any compact set $K\subset D$, there exists a positive constant
 $L_K>0$ such that:
$$ \parallel f(t,x)-f(t,y) \parallel \leq L_K \parallel x-y \parallel~,
~{\rm for~any~}x,y\in K~{\rm and}~t\in\R $$

4) bounded in time, that means there exists a continuous function
 $M:D\rightarrow\R$ such that:

$$ \parallel f(t,x) \parallel \leq M(x) ~~,~~ {\rm for~any}~ t\in\R$$

With these hypotheses we know that for any $(t_0,x_0)\in \R\times D$
there exists a unique solution of the Cauchy problem:
\begin{equation}
\label{e2}
\begin{array}{c}
\mbox{$\dot{x}=f(t,x)$} \\
\mbox{$x(t_0)=x_0$}
\end{array}
\end{equation}
with the initial data $(t_0,x_0)$; we denote by $x(t;t_0,x_0)$ this
solution. One can define this solution for $t\in(t_0-T,t_0+T)$ where 
$T=\sup_{r>0,B_r(x_0)\subset D}\frac{r}{\parallel f
{\parallel}_{B_r(x_0)}}$, the supremum is taken over all positive
 radius such that the ball centered around $x_0$, $B_r(x_0)=\{x\in
{\R}^n|\parallel x-x_0\parallel < r \}$, is completely included
 in $D$ and $\parallel f {\parallel}_{B_r(x_0)}=\sup_{(t,x)\in\R\times
 {\bar{B}}_r(x_0)}\parallel f(t,x) \parallel $ is a supremum norm of $f$ with
respect to $B_r(x_0)$ (where is no confusion we denote $B_r=B_r(0)$).
 The function ${\gamma}_{t,t_0}(x_0)=x(t;t_0,x_0)$
is well defined for some bounded open set $S$, ${\gamma}_{t,t_0}:S
\rightarrow U\subset D$ (with $U$ open and bounded) and it is Lipschitz
 continuous with a Lipschitz constant given by $L=exp(L_U |t-t_0|)$
 ($L_U$ being the Lipschitz constant associated to $f$, as above, on
 the compact set $\bar{U}$). All these results can be found in any
 textbook of  differential equations (for instance see \cite{cord}).

Our concern regards the stability behaviour of the equilibrium point
 $\bar{x}=0$. First we recall some definitions about stability (in
 Liapunov sense).

{\bf Definition} We say the equilibrium point $\bar{x}=0$ for (\ref{e1})
 is {\em uniformly stable}, if for any $\varepsilon >0$ there exists
 ${\delta}_{\varepsilon}>0$ such that for any $t_0\in\R$ and $x_0\in\R$
 with $\parallel x_0\parallel <{\delta}_{\varepsilon}$ the solution
 $x(t;t_0,x_0)$ is defined for all $t\geq t_0$ and furthermore
 $\parallel x(t;t_0,x_0)\parallel <\varepsilon$, for every $t>t_0$.

{\bf Definition} We say that the equilibrium point $\bar{x}=0$
 for (\ref{e1}) is {\em uniformly asymptotic stable}, if it is uniformly
 stable and there exists a ${\delta} >0$ such that for any  $t_0\in\R$ and
 $x_0\in D$ with $\parallel x_0\parallel < {\delta}$ the
 solution $x(t;t_0,x_0)$ is defined for every $t\geq t_0$ and $\lim_{t\rightarrow
\infty} x(t;t_0,x_0)=0$.

If in the  definition of uniform stability we interchange 'there exists 
${\delta}_{\varepsilon} >0$'  with 'for any $t_0\in\R$ ' (thus $\delta$
 will depend on $\varepsilon$ and $t_0$, ${\delta}_{\varepsilon,t_0}$)
 then the equilibrium is said (just) {\em stable}. If we proceed the same
 in the second definition we obtain that the equilibrium is {\em asymptotic
 stable}. For time-invariant systems there is no distinction between uniform
 stability and stability, or uniform asymptotic stability and asymptotic
 stability. In general case, the uniform (asymptotic) stability implies
 (asymptotic) stability, but the converse is not true (see for instance 
 \cite{b5}).

We say that the dynamics (\ref{e1}) has a {\em positive invariant set} $N$
 if for any $t_0\in \R$ and $x_0\in N$ the solution $x(t;t_0,x_0)\in N$
 for all $t\geq t_0$ for which it is well-defined. Then it makes sense to consider
 the {\em dynamics restricted to $N$}, i.e. the function:
 $$ X:{\R}^{+}\times\R\times N \rightarrow N~~~,~~~ X(\tau;t_0,x_0)=x(\tau+t_0;
t_0,x_0) $$
where $\tau$ runs up to a maximal value depending on $(t_0,x_0)$. 
Moreover, by considering the case of $f$ from (\ref{e1}) we obtain that
 $X(\tau;t_0,0)=0$, for any $\tau >0$, $t_0\in\R$. Therefore we may define
 the corresponding stability properties of the restricted dynamics
 as above, where we replace $D$ by $N$.

The main result of this paper is given by the following theorem:
\begin{Theorem}
\label{t1}

Consider the time-varying dynamical system $(\ref{e1})$ for which $f$ has
 the properties $1)-4)$. Suppose there exists a function $V:D\rightarrow\R$
 of class ${\cal C}^1$ such that:

$H1)$ $V(x)\geq 0$ for every $x\in D$ and $V(0)=0$;

$H2)$ There exists a continuous function $W:D\rightarrow\R$ such that:
$$\frac{dV}{dt}(t,x)=\nabla V (x)\cdot f(t,x) \leq W(x) \leq 0$$

$H3)$ Let $E=\{ x\in D |W(x)=0 \}$ denote the zero-set (or kernel) of $W$; 
 suppose that $f$ resticted to $E$ is time-invariant (i.e. $f(t,x)=f(t_0,x)$, 
 for every $t\in\R$
 and $x\in E$). Let us denote by $N$ the maximal positive invariant set in $E$
, i.e. for any $x_0\in N$ and $t_0\in\R$, $x(t;t_0,x_0)\in N$,
 for every $t\in [t_0,t_0+T_{x_0})$ in the maximal 
 interval of definition of the solution.

Then the dynamics $(\ref{e1})$ has at $\bar{x}=0$ an uniformly asymptotic
 stable equilibrium point if and only if
 the dynamics restricted to $N$ has an asymptotic stable equilibrium
 at $\bar{x}=0$.  $\Diamond$

\end{Theorem}

Even if it has appeared in the literature in a more general setting 
(I refer to \cite{yoshi}), it is worth mentioning the form the Invariance 
 Principle takes in this context:

\begin{Theorem}[Invariance Principle]
\label{tinv}

Consider the time-varying dynamical system $(\ref{e1})$ for which $f$ has
 the properties $1)-4)$. Suppose there exists a function $V:D\rightarrow\R$
 of class ${\cal C}^1$ such that:

$H1)$ It is bounded below, i.e. $V(x)\geq V_0$ for any $x\in D$ for some 
$V_0\in\R$;

$H2)$ There exists a continuous function $W:D\rightarrow\R$ such that:
$$\frac{dV}{dt}(t,x)=\nabla V (x)\cdot f(t,x) \leq W(x) \leq 0$$

$H3)$ Let $E=\{ x\in D |W(x)=0 \}$ denote the zero-set (or kernel) of $W$; 
 suppose that $f$
 resticted to $E$ is time-invariant (i.e. $f(t,x)=f(t_0,x)$, for any $t\in\R$
 and $x\in E$). Let us denote by $N$ the maximal positive invariant set 
 included in $E$, i.e. for any $x_0\in N$ and $t_0\in\R$, $x(t;t_0,x_0)\in N$,
 for any $t\in [t_0,t_0+T_{x_0})$ in the maximal 
 interval of definition of the solution.

Then any bounded trajectory of $(\ref{e1})$ tends to $N$, i.e. if $(t_0,x_0)$
 is the initial data for a bounded solution included in $D$~ then:
\begin{equation}
\label{eci}
\lim_{t\rightarrow\infty} d(x(t;t_0,x_0),N) = 0
\end{equation}

\end{Theorem}

{\bf Remarks} 1) There are two directions in which Theorem \ref{t1} 
 generalizes the well-known
 Barbashin-Krasovskii-LaSalle's Theorem (see \cite{lasal1}, \cite{lasal2}
 or \cite{b6}); firstly
 we require $V$ to be only nonnegative and not strictly positive, secondly
 we consider the time-varying dynamical systems. There exists in literature
 two earlier results in the first direction that I wish to comment.
 The first result that I am
 referring to is Lemma 5 from \cite{bro}. In this lemma only
 autononous systems are considered and the restricted dynamics is required to 
 be attractive in the sense that all trajectories should tend to the origin.
 I point out that only the requirement of attractivity is not enough;
 this can be seen in a trivial case, namely the 2 dimensional system
 given by Vinograd (conform \cite{b5}), for which the origin is an attractive
 equilibrium but not stable, and take $V\equiv 0$.
 I stress out that for the purposes
 of their paper (\cite{bro}) Lemma 5 can be replaced by Theorem \ref{t2} of this
 paper without affecting the other results of that paper.

A second result has appeared in \cite{b7} but not in a general and explicit form
 as here. In fact in \cite{b7} the author is concerned with the stability of the
 large-scale systems which are already decomposed in triangular form.
 Thus, this result
 solves the problem only in the case when we can perform the observability
 decomposition of the dynamics (\ref{e1}) with respect to the output $W(x)$.
 This case requires a supplementary condition, namely the codistribution span by
 $dW,dL_fW,\ldots,dL_f^nW$ to be of constant rank on $D$ (see \cite{b3}).
 Among other requirements, this geometric condition implies also that
 $N$ is a manifold, whereas we do not assume here this rather strong assumption.

 I acknowledge the existence of a recently published paper that deals with
 a similar extension of the Liapunov theorem, yet only for autonomous
  systems (\cite{pupu}). However, I was unaware of this result at the time
  I was working in this field (i.e. 1993--1995).

2) Some other papers deal with extensions of the invariance
 principle for nonautonomous systems. In two special cases, when the system
 is either asymptotically autonomous (in \cite{yoshi}) or 
asymptotically almost periodic (in \cite{miler}), the bounded solution tends to 
 the
 largest pseudo-invariant set in $E$. However
 they use the classical Liapunov theorems to obtain the uniform boundedness
 of the solutions. Thus they require the existence of a strictly positive 
 definite
 function playing the r\^{o}le of Liapunov function, while here we require
 only nonnegativeness of the Liapunov-like function. In other approaches
 an additional auxiliary function is assumed and by means of extra conditions
 the time in $E$ is controlled (see the results of Salvadori or Matrosov,
 e.g. in \cite{rouche}). In a third approach an extra condition on $\dot{V}$
 is considered without any additional condition on the vector field; such
 an approach is considered in \cite{aeyels}.

3) The condition that the restricted dynamics to be
 uniformly asymptoticaly stable is necessary and sufficient .
 Thus it is a center-manifold-type result
 where a knowledge about a restricted dynamics to some invariant set
 implies the same property of the whole dynamics. We point out here that the
set $N$ does not need to be a manifold.

4) One could expect that simple stability of the restricted dynamics
 would imply uniform stability of the restricted dynamics. But this is not
 true as we can see from the following example:
\begin{Example}
\label{exem1}
Consider the following autonomous planar system:
\begin{equation}
\label{ex}
\left\{  \begin{array}[c]{c}
\mbox{$\dot{x}=y^2$} \\
\mbox{$\dot{y}=-y^3$}
    \end{array} ~~,~~(x,y)\in {\R}^2
\right.
\end{equation}
The solution of the system is given by $(x,y)\rightarrow (x+ln(1+y^2 t),
\frac{y}{\sqrt{1+y^2 t}})$. It is obvious that the equilibrium is not stable
 but if we take $V=y^2$ we have $\frac{dV}{dt}=-2 y^4$ and on the set $E=N=\{
 (x,0),~x\in\R\}$  the dynamics is trivial stable $\dot{x}=0$. 
\end{Example}
The problem is not the nonisolation of the equilibrium, but the existence
 of some invariant sets in any neighborhood of the equilibrium;

5) Theorem \ref{tinv} is the natural generalization of the Invariance
 Principle to the class of systems considered in this paper. The conclusion
 of this theorem applies only to bounded trajectories. Thus we have to
 know apriori which solutions are bounded. Since they are bounded we can
 extend them indefinitely in positive time. Thus it makes sense to take
 the limit $t\rightarrow\infty$ in $(\ref{eci})$. We mention that a more general
 Invariance Principle  can be obtained even under
  weaker conditions than those from here (see \cite{yoshi}).

\vspace{4mm}

The organization of the paper is the following: in the next section
 we give the proof of these results. In the third section we consider the
 autonomous case and we present the systemic consequences related to the 
  nonlinear Liapunov equation  and  a special type of 
 zero-state detectability. 
 In the fourth section  we consider a nonlinear
 Riccati equation (or Hamilton-Jacobi equation) and we present a result
 of robust stabilizability by output feedback.
 The last section contains the conclusions and is followed by the bibliography.

\section{Proof of the Main Results}

We prove by contradiction  the uniform stability of the 
 equilibrium. For this, we construct a ${\cal C}^1$-convergent sequence of 
 solutions that are going away from the origin and whose limit is a 
 trajectory, thus contradicting the hypothesis.

For  the uniform asymptotic stability, we prove first that the $\omega$-limit
 set of bounded trajectories is included in $N$ (implicitely proving the
 Invariance Principle - Theorem \ref{tinv}) and then we addapt a classical
 trick (used for instance in Theorem 34.2 from \cite{b5}) that the
 convergence of trajectories in $\omega$-limit set will attract the
 convergence of the bounded trajectory itself. In both steps we use 
  essentially the time-invariant property of $f$ restricted to $E$.
 In proving the uniform stability we also obtain that the solution can 
 be defined
 on the whole positive real set (can be completely extended in future). 
 
Theorem \ref{tinv} (the Invariance Principle) will follow simply from a
 lemma that we state during the proof of uniform attractivity.
 
 First we need a lemma:

\begin{Lemma}
\label{l1}
 Let $f$ be a vector field defined on a domain $D$ and having the properties
 1-4 as above.
 Let $(t_i{)}_i$ be a sequence of real numbers and $(w_i{)}_i$, $w_i:[a,b]
\rightarrow D$ be a sequence of trajectories for the time-translated vector
 field $f$ with $t_i$, i.e. ${\dot{w}}_i(t)=f(t+t_i,w_i(t))$ . 

If the trajectories are uniformly bounded, i.e. there exists $M>0$ such
 that $\parallel w_i {\parallel}_{\infty}<M$, for any $i$,
 then we can extract a subsequence, denoted also by $(w_i{)}_i$, uniformly
 convergent to a function $w$ in ${\cal C}^1([a,b];D)$, i.e. $w_i\rightarrow w$
 and ${\dot{w}}_i\rightarrow \dot{w}$ both uniformly in ${\cal C}^0([a,b];D)$

\rm
{\bf Proof} We apply the Ascoli-Arzel\`{a} Lemma twice: first to extract a subsequence
 such that $(w_i{)}_i$ is uniformly convergent and second to extract further 
 another subsequence such that $({\dot{w}}_i{)}_i$ is uniformly convergent.
 Then we obtain that $\lim_i \frac{d}{dt} w_i = \frac{d}{dt} \lim_i w_i$.

1. We verify that $(w_i{)}_i$ are uniformly bounded and equicontinuous.
 The uniformly boundedness comes from $\parallel w_i {\parallel}_{\infty}<M$.
 The equicontinuity comes from the uniformly boundedness of the first
 derivative. Indeed, since $\parallel w_i \parallel \leq M$, the closed ball 
 ${\bar{B}}_M$
 is compact and $f(t,\cdot)$ is  continuous on ${\bar{B}}_M$, there
 exists a constant $A$ such that $\parallel f(t,x) \parallel \leq A$, for any
 $(t,x)\in\R\times {\bar{B}}_M$. Then:
$$\parallel {\dot{w}}_i(t) \parallel = \parallel f(t+t_i,w_i(t)) \parallel
 \leq A ~~,~{\rm for~any}~i~{\rm and}~t\in [a,b] $$
Thus $(w_i{)}_i$ is relatively compact and we can extract a subsequence,
 that we denote also by $(w_i{)}_i$, which is uniformly convergent to
 a function $w\in {\cal C}^0([a,b];D)$.

2. We prove that $({\dot{w}}_i{)}_i$ is relatively compact.
 We have already proved the uniform boundedness $\parallel {\dot{w}}_i 
{\parallel}_{\infty}
 \leq A$. For the equicontinuity we use both the uniform continuity in $t$
 and uniform local Lipschitz continuity in $x$, of $f$. Let $L_M$ be the
 uniform Lipschitz constant corresponding to the compact set ${\bar{B}}_M$.
 Then:
$$\parallel {\dot{w}}_i(t_1) - {\dot{w}}_i(t_2) \parallel = \parallel
 f(t_i+t_1,w_i(t_1))- f(t_i+t_2,w_i(t_2))\parallel \leq $$
$$\parallel f(t_i+t_1,w_i(t_1))-f(t_i+t_2,w_i(t_1)) \parallel$$ 
$$+ \parallel
 f(t_i+t_2,w_i(t_1))-f(t_i+t_2,w_i(t_2)) \parallel $$
Let $\varepsilon >0$ be arbitrarily. Then we choose ${\delta}_1$ such that
$ \parallel f(s_1,x)-f(s_2,x)\parallel <\frac{\varepsilon}{2}~~,
~~{\rm for~any}~|s_1-s_2|<{\delta}_1~{\rm and}~x\in{\bar{B}}_M $
On the other hand: $\parallel f(t_i+t_2,w_i(t_1))-t(t_i+t_2,w_i(t_2))\parallel 
\leq
 L_M \parallel w_i(t_1) - w_i(t_2) \parallel \leq L_M A |t_1-t_2 |$. 
Then we choose  $\delta=\min({\delta}_1,\frac{\varepsilon}{2 L_M A})$. Then
 the left-hand side from the above inequality is also bounded by 
$\frac{\varepsilon}{2}$ for any $t_1,t_2$ with $|t_1-t_2|<\delta$. Thus
 $\parallel {\dot{w}}_i(t_1)-{\dot{w}}_i(t_2) \parallel < \frac{\varepsilon}{2}
 + \frac{\varepsilon}{2} = \varepsilon $, for any $i$ and $t_1,t_2\in [a,b]$,
 $|t_1-t_2|<\delta$.

We can now extract a second subsequence from $(w_i{)}_i$ such that $({\dot{w}}_i
{)}_i$ is also uniformly convergent and this ends the proof of lemma.  $\Box$
\end{Lemma}

\vspace{4mm}

{\bf Proof of Uniform Stability} 

Let us assume that the equilibrium is not uniformly stable. Then there exists
 ${\varepsilon}_0>0$ such that for any $\delta$, $0<\delta<{\varepsilon}_0$
 there are $x_0,t$ and $\Delta>0$ such that $\parallel x_0\parallel <\delta$
 and $\parallel x(t+\Delta;t,x_0)\parallel = {\varepsilon}_0$, $\parallel
 x(t+\tau;t,x_0)\parallel <{\varepsilon}_0$, for $0\leq \tau <\Delta$.
 We choose ${\varepsilon}_0$ (eventually by shrinking it) such that
 $\bez\cap N$ is included in the attraction domain
 of the origin (for the restricted dynamics).

By choosing a sequence $({\delta}_i{)}_i$ converging to zero we obtain
 sequences $(x_{0i}{)}_i$, $(t_i{)}_i$ and $({\Delta}_i{)}_i$ such that:
$\parallel x_{0i} \parallel \rightarrow 0$  and
$\parallel x(t_i+{\Delta}_i;t_i,x_{0i}) \parallel ={\varepsilon}_0$

Let ${\delta}<{\varepsilon}_0$ be such that for any $z_0\in B_{\delta}\cap N$ 
we have \linebreak 
$\parallel x(t;0,z_0)\parallel <\frac{{\varepsilon}_0}{2}$ for any $t>0$
 (such
 a choice for $\delta$ is possible since the dynamics restricted to $N$ is
 stable). Let $i_0$ be such that ${\delta}_i<\delta$, for $i>i_0$. 
 We denote by $(u_i)_{i>i_0}$ the time moments such that $\parallel 
 x(t_i+u_i;t_i,
 x_{0i})\parallel =\delta$ and $\parallel x(t;t_i,x_{0i})\parallel >\delta$
 for $t>t_i+u_i$.
 Since the spheres ${\bar{S}}_{{\varepsilon}_0}$ and ${\bar{S}}_{\delta}$
 are compact we can extract
 a subsequence (indexed also by $i$) such that both
 $x_i=x(t_i+{\Delta}_i;t_i,x_{0i})$ and $y_i=x(t_i+u_i;t_i,x_{0i})$ are
 convergent to $x^*$, respectively to $y^*$; $x_i\rightarrow x^*$,
 $y_i\rightarrow y^*$, $\parallel x^*\parallel ={\varepsilon}_0$,
 $\parallel y^*\parallel =\delta$. Since $V$ is
 continuously nonincreasing on trajectories and $\lim_i V(x_{0i})=0$, we get 
 $V(x^*)=V(y^*)=0$. Therfore $x^*,y^*\in N$.
 
 Suppose $\parallel f(x,t)\parallel \leq A$ on $\bez$,
 for some $A>0$. Then one can easily prove  that ${\Delta}_i-u_i\geq \frac{
{\varepsilon}_0-\delta}{A}=T_1$, for any $i>i_0$ (i.e. the flight time between
 two spheres of radius $\delta$ and ${\varepsilon}_0$ has a lower bound).

Define now the time-translated vector fields $f_i(t,x)=f(t+t_i+u_i,x)$ and
 denote by $w_i:[0,T_1]\rightarrow \bez$ the
 time-translated solutions $w_i(t)=x(t+t_i+u_i;t_i,x_{0i})$. Then:
$ {\dot{w}}_i(t)=f_i(t,w_i(t))$, $0\leq t\leq T_1$.
By applying Lemma \ref{l1} we get a subsequence uniformly convergent to a
 trajectory $w^1:[0,T_1]\rightarrow \bez\cap N$,
 such that $w^1(0)=\lim_iw_i(0)=y^*$ and $\parallel w^1(t)\parallel >\delta$,
 for $0<t\leq T_1$. If $\parallel w^1(T_1)\parallel <{\varepsilon}_0$ we
 obtain that ${\Delta}_i-u_i-T_1 > \frac{{\varepsilon}_0-\parallel w^1(T_1)
 \parallel }{A}$, for some $i\geq i_1>i_0$.  Then, we denote $T_2=T_1+
 \frac{{\varepsilon}_0-\parallel w^1(T_1)\parallel }{A}$ and we repeat the
 scheme. We obtain another sequence which is uniformly convergent to a
 trajectory $w^2:[0,T_2]\rightarrow \bez\cap N$
 such that $w^2(0)=y^*$, $\parallel w^2(t)\parallel >\delta$, $0< t\leq T_2$
 and $w^2(t)=w^1(t)$, for $0\leq t\leq T_1$.

Thus we extend each trajectory $w^k:[0,T_k]\rightarrow \bez\cap N$
 to a trajectory $w^{k+1}:[0,T_{k+1}]\rightarrow \bez\cap N$ such that
 $T_{k+1}\geq T_k$, $w^{k+1}(t)=w^k(t)$ for $0\leq t\leq T_k$ and
 $\parallel w^{k+1}(t)\parallel >\delta$, for $0 < t\leq T_{k+1}$.

We end this sequence of extensions in two cases:

1) $\lim_k T_k = T^*<+\infty$ (the limit may be reached in a finite number
 of steps), in which case we have\linebreak 
 $\lim_k\parallel w^k(T_k)\parallel={\varepsilon}_0$ and thus 
 $\lim_kw^k(T_k)=x^*$; or:

2) $\lim_kT_k=+\infty$.

In the first case we obtain a trajectory $w^*:[0,T^*]\rightarrow \bez\cap N$
 such that $w^*(0)=y^*$, $w^*(T^*)=x^*$ with $\parallel w^*(0)\parallel =\delta$
 and $\parallel w^*(T^*)\parallel ={\varepsilon}_0$. But this is a contradiction
 with the choice of $\delta$ (and of stability of the restricted dynamics).

In the second case we obtain a trajectory $w^*:[0,\infty)\rightarrow\bez\cap N$
 such that $\parallel w^*(0)\parallel =\delta <{\varepsilon}_0$ and
 $\parallel w^*(t)\parallel >\delta$ for $t>0$. Thus $\lim_{t\rightarrow\infty}
 w^*(t)\neq 0$ contradicting the assumption that $\bez\cap N$ is included in the
 attraction domain of the origin. Now the proof is  complete.  $\Box$.
\vspace{4mm}

For the proof of uniformly attractivity we recall a few definitions and results:

{\bf Definition} A point $x^*$ is called {\em $\omega$-limit point} for
 the trajectrory $x(t;t_0,x_0)$ if there exists a sequence $(t_k{)}_k$ such that
 $\lim_{k\rightarrow\infty} t_k=\infty$, $x(t;t_0,x_0)$ is defined for all $t>t_0$
 and $\lim_k x(t_k;t_0,x_0)=x^*$.

The set of all $\omega$-limit points is called the $\omega$-limit set and
 is denoted by $\Omega(t_0,x_0)$. It characterizes the trajectory $x(t;t_0,x_0)$
 and it depends on the initial data $(t_0,x_0)$.

\begin{Theorem}[Birkoff's Limit Set Theorem, see \cite{bir}]
\label{tbir} 
 A bounded trajectory approaches its
 $\omega$-limit set, i.e. \linebreak 
 $\lim_{t\rightarrow\infty} d(x(t;t_0,x_0),\Omega(t_0,x_0))=
 0$, where $d(p,S)=\inf_{x\in S}\parallel p-x \parallel $ is the distance
 between the point $p$ and the set $S$. $\Diamond$
\end{Theorem}

There is also a very useful result about uniformly continuous functions:

\begin{Lemma}[Barb\u{a}lat's Lemma, see \cite{b4}]
 If $g:[t_0,\infty)\rightarrow\infty$ is a
 uniformly continuous function such that the following limit exists and is
 finite, $\lim_{t\rightarrow\infty}\int_{t_0}^{t} g(\tau)d\tau$, then
 $\lim_{t\rightarrow\infty} g(t)=0$.
 $\Diamond$
\end{Lemma}

\vspace{4mm}

{\bf Proof of Uniform Attractivity}

We already know that $\bar{x}=0$ is uniformly stable. What we have to prove is
 the uniform attractivity. 

Let ${\varepsilon}_0>0$ be chosen with the following properties: for any
 $t_0$ and $x_0\in D\cap\bez$ the positive trajectory $x(t;t_0,x_0)$ is bounded
 by ${\varepsilon}_1$ (i.e. $x(t;t_0,x_0)\in B_{{\varepsilon}_1}$); for any
 $t_1$ and $x_1\in D\cap B_{{\varepsilon}_1}$ the trajectory
 $x(t;t_1,x_1)$,  $t>t_1$, is bounded by some $M$; and for any 
 $x_2\in N\cap B_{{\varepsilon}_1}$ the trajectory $x(t;t_0,x_2)$ tends to the
 origin $\lim_{t\rightarrow\infty} x(t;t_0,x_2)=0$. We are going to prove that
 $\lim_{t\rightarrow\infty} x(t;t_0,x_0)=0$.

Let us consider the $\omega$-limit set $\Omega (t_0,x_0)$. It is enough to prove
 that $\Omega (t_0,x_0)=\{ 0 \}$, because of Birkoff's Limit Set Theorem.

Let $x^*\in\Omega(t_0,x_0)$ and suppose $x^*\neq 0$. Let us denote by
$x(t)=x(t;t_0,x_0)$ and $g(t)=\nabla V(x(t))\cdot f(t,x(t))$. Since the
 solution is continuous and bounded, so is $g(t)$. On the other hand
 $$ V(x(t))=V(x_0)+\int_{t_0}^t g(\tau)d\tau $$
Since $\dot{x}(t)=f(t,x(t))$ and $x(t)$ is bounded we obtain that it is also
 uniformly continuous. Thus $g(t)$ is also uniformly continuous (recall we have
 assumed $f(\cdot,x)$ is uniformly continuous in $t$). Let $(t_k{)}_k$ be a
 sequence that renders $x^*$ a $\omega$-limit point. Then $\lim_k V(x(t_k))=
 V(\lim_k x(t_k))=V(x^*)$. Since $V(x(t))$ is a decreasing function bounded
 below, there exists the limit: $\lim_{t\rightarrow\infty}V(x(t))=V(x^*)$.
 Now, applying Barb\u{a}lat's Lemma we obtain $\lim_{t\rightarrow\infty}g(t)=0$
 or $W(x^*)=0$. Thus $\Omega(t_0,x_0)\subset E$, the kernel of $W$.

In this point we need a result about the behaviour of solutions starting at
 $x^*$. We mention that the following lemma is a consequence of Theorem 3
 from \cite{yoshi}. But, since we are under stronger conditions, we have
 found a simpler proof that we are going to present here (our conditions
 are stronger because we need to obtain uniform stability and consequently
 boundedness of the solutions when Liapunov function is only positive
 semidefinite, which overall means a weaker condition). 

\begin{Lemma}
\label{l2}
The positive trajectory starting at $x^*$ is included in $E$ and thus the
 $\Omega$-limit set is a positive invariant set included in $N$.

\rm
{\bf Proof}

Let $\tau>0$ be an arbitrary time interval. Let $(t_k)_k$ be the sequence
 that renders $x^*$ a $\omega$-limit point for the trajectory $x(t)=x(t;
t_0,x_0)$. Then, if we denote by $x_k=x(t_k)$ we have $\lim_k x_k=x^*$. Consider
 the following sequence of functions:
$w_k:[0,\tau]\rightarrow D~~,~~ w_k(t)=x(t+t_k;t_k,x^*)$. 
 We have chosen $x_0,\,t_0$ such that all these functions are bounded by $M$, i.e.
 $\parallel w_k {\parallel }_{\infty}<M$. We have $w_k(0)=x^*$ and $V(w_k(t))\leq
 V(x^*)$. Let us denote by $y_k^t=x(t+t_k)$, for any $0\leq t\leq \tau$, and
 let $L_M$ be the Lipschitz constant of $f$ on the compact ${\bar{B}}_M$. Then:
 $\parallel y_k^t - w_k(t) \parallel \leq e^{L_m t}\parallel x_k - x^* 
 \parallel$
and, since $\lim_k x_k=x^*$ we get $\lim_k\parallel y_k^t -w_k(t) \parallel =0$.
 On a hand, since $V(x^*)=\lim_{t\rightarrow\infty}V(x(t))$ and $V$ is 
nonincreasing on trajectories we have $V(y_k^t)>V(x^*)$ and also $\lim_k
 V(y_k^t)=V(x^*)=\lim_k V(w_k(t))$. On the other hand, since $(w_k)_k$ are 
uniformly bounded we apply Lemma \ref{l1} and we obtain a subsequence uniformly
 convergent to a function $w\in {\cal C}^1([0,\tau];D\cup {\bar{B}}_M)$. 
 Obviously $V(w(t))=V(x^*)$ for any $0\leq t\leq\tau$. Thus $W(w(t))=0$ and
 $w(t)\in E$. On the other hand, since $f$ is continuous in $(t,x)$ we obtain
 that $w$ is an integral curve of $f$, i.e. $\dot{w}(t)=f(t_*,w(t))$, for
 $0\leq t\leq\tau$ and any $t_*$. In particular, for $t_*=t_k$ we get $w(t)$
 is a solution of the same equation as $w_k(t)$ and $w(0)=w_k(0)=x^*$. By the
 uniqueness of the solution they must coincide. Then $x(t+t_k;t_k,x^*)\in E$
 for $0\leq t\leq\tau$. But $\tau$ was arbitrarily; thus $x(t;t_0,x^*)\in E$ for
 any $t$ and then $x^*\in N$.   $\Box$
\end{Lemma}

Since the trajectory starting at $x^*$ is included in $N$, it should converge
 to the origin (the equilibrium point). Let us denote by $\varepsilon = \frac{
\parallel x^* \parallel }{2}$. From uniform stability there exists a $\delta>0$
 such that for any $\tilde{x}\in D$, $\parallel\tilde{x}\parallel <\delta$
 implies $\parallel x(t_2;t_1,\tilde{x}) \parallel <\varepsilon$, for any
 $t_2>t_1$. Let $\Delta t$ be a time interval
 such that $\parallel x(t;0,x^*)\parallel < \frac{\delta}{2}$ for any $t>\Delta t$.
 We consider the compact set $C$, the $\frac{\delta}{2}$-neighborhood of the
 compact curve $\Gamma=\{ x(t;0,x^*)|0\leq t\leq \Delta t\}$:
$$ C=\{ x\in D | d(x,\Gamma)\leq \frac{\delta}{2} \} = \bigcup_{t\in [0,\Delta t]}
 \bar{B_{\delta /2}(x(t;0,x^*))} $$
which is the union of the closed balls centered at $x(t;0,x^*)$ and of radius
 $\frac{\delta}{2}$. We set ${\delta}_1=\frac{\delta}{2} exp(-L_C \Delta t)$
 where $L_C$ is the uniform Lipschitz constant of $f$ on the compact set $C$.
 Since the solution is uniformly Lipschitz with respect to the initial point
 $x_0$ we have that for any $t_1\in\R$ and $x_1$ such that $\parallel x_1 - x^*
 \parallel < {\delta}_1$ we get:
$\parallel x(t_1+\Delta t;t_1,x_1) - x(\Delta t;0,x^*) \parallel < 
\frac{\delta}{2}$
 and then $\parallel x(t_1+\Delta t;t_1,x_1) \parallel < \delta$. Furthermore,
 from the choice of $\delta$ we obtain that $\parallel x(t_1+\tau;t_1,x_1)
 \parallel <\varepsilon$, for any $\tau>\Delta t$ or
 $\parallel x(t_1+\tau;t_1,x_1) - x^* \parallel >\varepsilon$, for any 
 $\tau>\Delta t$.

Now we pick a $t_n$ such that $\parallel x(t_n;t_0,x_0) -x^* \parallel < 
{\delta}_1$. Then, from the previous discussion $\parallel x(t_n+\tau;t_0,x_0)-
 x^*\parallel >\varepsilon$, for any $\tau>\Delta t$ which contradicts the limit
 $\lim_k x(t_k;t_0,x_0)=x^*$. This contradiction comes from the hypothesis
 that $x^*\neq 0$. Thus $\Omega(t_0,x_0)=\{ 0 \}$ and now the proof is complete.
 $\Box$
\vspace{4mm}

{\bf Proof of Theorem \ref{tinv} (The Invariance Principle)}

If $x(t;t_0,x_0)$ is a bounded trajectory then, from Birkoff's Limit Set
 Theorem it approaches its $\omega$-limit set. On the one hand we can use
 Barb\u{a}lat's Lemma and prove that $W$ vanishes on $\omega$-limit set
 of bounded trajectories. On the other hand, as we have proved in Lemma
 \ref{l2}, the $\omega$-limit set is invariant and included in $N$. Thus
 the bounded trajectory approaches the set $N$. $\Box$

\section{The Autonomous Case: Consequences in Nonlinear Control Theory}

 Consider the following inputless nonlinear control system:

\begin{equation}
\label{e31}
S \left\{ \begin{array}[c]{c}
\mbox{$\dot{x}=f(x)$} \\
\mbox{$y=h(x)$} 
     \end{array} \right. ~~,~x\in D\subset{\R}^n~y\in{\R}^p
\end{equation}
such that $f(0)=0$, $h(0)=0$ and $D$ a neighborhood of the origin. Suppose $f$
 is local Lipschitz continuous and $h$ continuous on $D$. Then denote by
 $x(t,x_0)$ the flow generated by $f$ on $D$ (i.e. the solution of $\dot{x}=
 f(x)$, $x(0)=x_0$), by $E=ker~h=\{x\in D|h(x)=0\}$, the kernel of $h$
 and by $N$ the maximal positive invariant set included in $E$, i.e.
 the set $N=\{ \tilde{x}\in D |h(x(t,\tilde{x}))=0~{\rm for~any}~t\geq 0~{\rm
 such~that}~x(t,\tilde{x})~{\rm has~sense}\}$.

We present two concepts of detectability for $(\ref{e31})$.
 The first one has been used by many authors lately 
(see for instance \cite{isi}).

{\bf Definition} The pair $(h,f)$ is called {\em zero-state detectable}
 (or {\em z.s.d.}) 
 if $\bar{x}=0$ is an attractive point for the dynamics restricted to $N$,
 i.e. there exists an ${\varepsilon}_0>0$ such that for any $x_0\in N$,
 $\parallel x_0 \parallel <{\varepsilon}_0$, 
 $\lim_{t\rightarrow\infty} x(t,x_0)=0$.

{\bf Definition} The pair $(h,f)$ is called {\em strong zero-state detectable}
 (or {\em strong z.s.d.}) 
 if $\bar{x}=0$ is an asymptotical stable equilibrium  point for the dynamics
 restricted to $N$,
 i.e. it is zero-state detectable and for some ${\varepsilon}_0$ and
 for any $x_0\in N$ with $\parallel x_0\parallel < {\varepsilon}_0$,
 $\lim_{t\rightarrow\infty} x(t,x_0)=0$.

We see that strong z.s.d. implies z.s.d., but obviously 
the converse is not true.
 
In this framework, as a consequence of the main result we can state the following theorem:

\begin{Theorem}
\label{t2}
For the inputless nonlinear control system $(\ref{e31})$ with $f$ local Lipschitz continuous and $h$ continuous, consider the following
 nonlinear Liapunov equation:
\begin{equation}
\label{e32}
\nabla V \cdot f +\parallel h {\parallel }^q =0
\end{equation}
or the following nonlinear Liapunov inequality:
\begin{equation}
\label{e33}
\nabla V \cdot f +\parallel h {\parallel }^q \leq 0
\end{equation}
for some $q>0$. Suppose there exists a positive semidefinite solution of $(\ref{e32})$ or
 $(\ref{e33})$ of class ${\cal C}^1$ defined on $D$ such that $V(0)=0$.

 Then  the pair $(h,f)$ is strong zero-state detectable if and only if
  $\bar{x}=0$ is an
 asymptoticaly stable equilibrium for the dynamics $(\ref{e31})$. $\Diamond$
\end{Theorem}

Below we give an example:
\begin{Example}
\label{ex2}

Consider the dynamics:
\begin{equation}
\label{eheh}
\begin{array}{l}
\mbox{${\dot{x}}_1 =-x_1^3+\Psi(x_2)$} \\
\mbox{${\dot{x}}_2=-x_2^3$}
\end{array} ~~,~(x_1,x_2)\in{\R}^2
\end{equation}
where $\Psi:\R\rightarrow\R$ is local Lipschitz continuous, $\Psi(0)=0$
 and there exist constants $a>0$, $b\geq 1$ such that:
 $$ |\Psi(x)|\leq a|x|^b~~,~\forall x_2 $$
 If we choose as output function $h(x)=x_2^2$ we see that the pair $(h,f)$
  is strong zero-state detectable; indeed, the set $E=\{ x\in{\R}^2|h(x)=0\}
   = \{ (x_1,0)|x_1\in\R \}$ and the dynamics restricted to $E$ is 
    ${\dot{x}}_1=-x_1^3$ which is asymptoticaly stable.

Now, if we choose $ V(x)=\frac{x_2^2}{2}$
 we have $\dot{V}=-x_2^4$ and thus $V$ is a solution of the Liapunov
 equation $(\ref{e32})$ with $q=2$. Then, the equilibrium is asymptoticaly
 stable, as a consequence of the theorem \ref{t2}.

On the other hand we can explicitely solve for $x_2$: 
$x_2(t)=\frac{x_{20}}{\sqrt{2(1+x_{20}^2t)}}$
and then we have: $|\Psi(x_2(t))|\leq C(1+Bt)^{-1/2}$ for some $B,C>0$ and any
 $t\geq 0$. Now the asymptotic stability follows as a consequence of Theorem
 68.2 from \cite{b5} (stability under perturbation).

\end{Example}

\section{An Application in Robust Stabilizability}

 We present here, as an application, a robust stabilizability result
 for a nonlinear affine control system. In fact it is an absolute stability
 result about a particular situation. More general results about absolute stability
 for nonlinear affine control system will appear in a forthcoming paper.
 We base our approach on the existance of a positive semidefinite solution
 of some Hamilton-Jacobi equation or inequality. 
 Discussions about solutions of this type of equation may be found
 in \cite{schaft}.

Consider the following Single Input - Single Output control system:
\begin{equation}
\label{affine}
\left\{ \begin{array}{rcl}
 \mbox{$\dot{x}$} & = & f(x)+g(x)u \\
   y & = & h(x)
         \end{array} \right. ~~,~~x\in D\subset{\R}^n ~,~u,y\in\R
\end{equation}
where $f$ and $g$ are local Lipschitz continuous vector fields on a domain
 $D$ including the origin, $h$ is a local Lipschitz real-valued function
 on $D$, and $f(0)=0$, $h(0)=0$. Consider also a local Lipschitz output
 feedback:
\begin{equation}
\label{feed}
\varphi:\R\rightarrow\R~~,~~\varphi(0)=0
\end{equation}
We define now two classes of perturbations associated to this feedback.
 Let $a>0$ be a positive real number. The first class contains time-invariant
 perturbations:
$$P_1=\{ p:\R\rightarrow\R ~,~p~{\rm is~local~Lipschitz},~p(0)=0~{\rm and}~
|p(y)|<a|\varphi(y)|~,~\forall y\neq 0\}$$
while the second class is composed by time-varying perturbations:
$$P_2=\{ p:\R\times \R\rightarrow\R~,~p(y,t)~{\rm is~local~Lipschitz~in}~y~{\rm
 for}~t~{\rm fixed~and~uniformly~continuous~in}~t$$
$${\rm for~any}~y~{\rm fixed}~,~p(0,t)\equiv 0
~{\rm and~there~exists}~\varepsilon>0~{\rm such~that}~|p(y,t)|<(a-\varepsilon)
 |\varphi(y)|~,~\forall y\neq 0,t\}$$

Now we can define more precisely the concept of robust stability:

{\bf Definition} We say  the feedback $(\ref{feed})$ {\em robustly stabilizes}
 the system $(\ref{affine})$ {\em with respect to the class $P_1\cup P_2$}
 if for any perturbation $p\in P_1\cup P_2$ the closed-loop with the 
 perturbed feedback $\varphi +p$ has an asymptoticaly stable equilibrium at the
 origin.

In other words, we require that the origin to be asymptoticaly stable
 for the dynamics:
\begin{equation}
\label{dyn}
\dot{x}=f(x)+g(x) (\varphi(h(x))+p(h(x),t))
\end{equation}
for any $p\in P_1\cup P_2$. Since the null function belongs to $P_1$, the
 feedback $\varphi$ itself must stabilize the closed-loop too.

With these preparations we can state the result:

\begin{Theorem}
\label{t7}
 Consider the nonlinear affine control system $(\ref{affine})$ and the
 feedback $(\ref{feed})$. Suppose the pair $(h,f)$ is strong zero-state
 detectable and suppose the following Hamilton-Jacobi equation:
\begin{equation}
\label{ricc}
\nabla V\cdot f + (\frac{1}{2}\nabla V\cdot g + \varphi\circ h)^2-
 (1-a^2)(\varphi\circ h)^2 =0 ~~,~~V(0)=0
\end{equation}
or inequality:
\begin{equation}
\label{ciuciu}
\nabla V\cdot f + (\frac{1}{2}\nabla V\cdot g + \varphi\circ h)^2-
 (1-a^2)(\varphi\circ h)^2 \leq 0 ~~,~~V(0)=0
\end{equation}
  has a positive semidefinite solution $V$ of class ${\cal C}^1$ on $D$.

Then the feedback $\varphi$ robustly stabilizes the system $(\ref{affine})$
 with respect to the class $P_1\cup P_2$.

{\bf Proof}

Let us consider a perturbation $p\in P_1\cup P_2$. Then, the closed-loop
 dynamics is given by $(\ref{dyn})$. We compute the time derivative of
 the solution $V$ of $(\ref{ricc})$ with respect to this dynamics:
$$\frac{dV}{dt}=\nabla V \cdot f(x) + \nabla V\cdot g(x) (\varphi(h(x))+
 p(h(x),t))$$
After a few algebraic manipulations we get:
$$\frac{dV}{dt}\leq -(\frac{1}{2}\nabla V \cdot g - p\circ h)^2 +
 (p\circ h)^2- a^2 (\varphi\circ h)^2 $$
Now, for $p\in P_1$, $\frac{dV}{dt}$ is time-independent and we may take
 for instance:
$$W(x)=(p(h(x)))^2-a^2(\varphi(h(x)))^2 \leq 0$$
For $p\in P_2$, $\frac{dV}{dt}$ is time-dependent and we define:
$$W(x)= - (2a\varepsilon - {\varepsilon}^2)(\varphi(h(x)))^2\leq 0$$
Either a case or the other, we obtain (recall the definitions of $P_1$
 and $P_2$):
$$ \frac{dV}{dt} \leq W(x) \leq 0 $$
The kernel-set of $W$ is given by:
$$E=\{ x\in D|~W(x)=0 \} = \{ x\in D | ~h(x)=0 \} $$
We see that the closed-loop dynamics $(\ref{dyn})$ restricted to $E$
 is simply given by $\dot{x}=f(x)$ and is time-independent. Moreover, since
 we have supposed $(h,f)$ is strong zero-state detectable, it follows that
 the restricted dynamics to the maximal positive invariant set in $E$
 has an asymptoticaly stable equilibrium at the origin. Now,
 applying Theorem \ref{t1}, the result follows.  $\Box$
\end{Theorem}
 
Let us consider now an example:

\begin{Example}
\label{exe3}

Consider the following planar nonlinear control system:
\begin{equation}
\label{exec1}
\left\{ \begin{array}{l}
 \mbox{${\dot{x}}_1=-x_1^3+u$} \\
 \mbox{${\dot{x}}_2=-x_2^3$} \\
 \mbox{$y=x_2^3$}
       \end{array}  \right.
\end{equation}
We are interested to find how robust the feedback $\varphi(y)=y$ is, i.e.
 how large we can choose $a$ such that $\varphi$ robustly stabilizes
  the system $(\ref{exec1})$ with respect to the class $P_1\cup P_2$.

The Hamilton-Jacobi equation $(\ref{ricc})$ takes the form:
$$-x_1^3 \frac{\partial V}{\partial x_1} - x_2^3\frac{\partial V}{\partial
 x_2} +(\frac{1}{2}\frac{\partial V}{\partial x_1}+x_2^3)^2 - (1-a^2)
  x_2^6 =0 $$
or:
$$-x_1^3\frac{\partial V}{\partial x_1} - x_2^3 \frac{\partial V}{\partial
 x_2} +\frac{1}{4} (\frac{\partial V}{\partial x_1})^2 + x_2\frac{\partial
  V}{\partial x_1}+a^2 x_2^6 =0 $$
A solution of this equation is:
$$ V(x_1,x_2)=\frac{a^2}{4} x_2^4 $$
For any $a>0$ it is positive semidefinite and the system $(\ref{exec1})$
 is strong zero-state detectable. Thus, as a consequence of theorem
  \ref{t7}, we can choose $a$ arbitrary large such that $\varphi$ robustly
 stabilizes the system $(\ref{exec1})$ with respect to the class $P_1\cup
 P_2$.

 On the other hand, for any feedback $\Phi$, local Lipschitz and:
 $$|\Phi(y)|\leq a |y| ~~,~{\rm for~some}~a>0 $$
 we have seen in the previous example that the closed-loop has an
 asymptoticaly stable equilibrium at the origin.
\end{Example}

\section{Conclusions}

In this paper we study an extension of Barbashin-Krasovski-LaSalle and
 Invariance Principle to a class of time-varying dynamical systems. We
 impose two type of conditions on the vector field: one is regularity
 (we require uniformly continuity with respect to $t$ and uniformly local
 Lipschitz continuoity and boundedness with respect to $x$); the other 
 condition  requires the vector field to be time-invariant on the zero-set
 $E$ of an auxiliary function. In this setting
 we find that the asymptotic behaviour of the dynamics restricted to the
 largest positive invariant set in $E$ determines the asymptotic stability
 character of the full dynamics. 

 Then we study two applications in control theory. The first application
 concerns the notion of detectability. We give  another definition
 for this notion, called strong zero-state detectability 
  and we show how the existence of a positive semidefinite
  solution of the Liapunov equation or inequation is related to the asymptotic 
 stability
  of the equilibrium. We obtain a nonlinear equivalent of the linear well-known
  result: if the pair $(C,A)$ is detectable and there exists a positive solution
  $P\geq 0$ of the Liapunov algebraic equation $A^TP+PA+C^TC=0$, then the
  matrix $A$ has all eigenvalues with negative real part.

The second application is on the problem of robust stabilizability. We give
 sufficient conditions such that a given feedback robustly stabilizes the
 closed-loop with respect to two sector classes of perturbations (time-invariant
 and time-varying). The condition is formulated in term of the existence
 of a positive solution of some Hamilton-Jacobi equation or inequality.

 This last application opens the problem of absolute stability for nonlinear
  affine control systems, that will be considered in a forthcoming paper.

\end{document}